\documentclass[12pt]{article}
\usepackage{amssymb}
\usepackage{amsthm}
\usepackage{amsmath}
\usepackage{graphicx}
\graphicspath{{./Figs/}}
\usepackage{lmodern}
\usepackage{graphics}
\usepackage{float}

\newtheorem{teo}{Theorem}
\newtheorem{lem}{Lemma}

\newcommand{\pa}{\partial}

\newcommand{\vp}{\varphi}
\newcommand{\ve}{\varepsilon}
\newcommand{\om}{\omega}
\newcommand{\be}{\begin{equation}}
\newcommand{\ee}{\end{equation}}

\newcommand{\ipd}{\stackrel{\normalfont\text{def}}{=}}
\newcommand{\const}{\operatorname{const}}
\newcommand{\dpsi}{\dot{\psi}_0}

\usepackage[usenames]{color}
\begin{document}
\allowdisplaybreaks

\title{Soliton dynamics for  the  general Degasperis-Procesi equation }

\author{ Georgy~Omel'yanov\thanks{
Universidad de Sonora, Rosales y Encinas s/n, 83000, Hermosillo,
Sonora, M\'exico,\ omel@mat.uson.mx} }
\date{}
\maketitle
\begin{abstract}
We consider  the general Degasperis-Procesi model of shallow water out-flows. This fife parametric family of conservation laws contains, in particular,   KdV, Camassa-Holm, and Degasperis-Procesi    equations.  The main result consists of a criterion which guarantees the existence of a smooth soliton type solution. We discuss also the scenario of soliton interaction for this model in the nonintegrable case.
\end{abstract}


\section{Introduction}
 The general Degasperis-Procesi model (\cite{ DegProc}, 1999) is the fife parametric family of conservation laws
\begin{align}
&\frac{\pa }{\pa t}\left\{u-\alpha^2\ve^2\frac{\pa^2 u}{\pa x^2}\right\}\label{1}\\
&+\frac{\pa}{\pa x}\left\{c_0u+c_1u^2
-c_2\ve^2\Big(\frac{\pa u}{\pa x}\Big)^2+\ve^2\big(\gamma-c_3u\big)\frac{\pa^2 u}{\pa x^2}\right\}=0, \; x \in \mathbb{R}^1, \; t > 0,\notag
\end{align}
which describes, in particular,  the dynamics of out-flows of shallow water.
Here $\alpha$, $c_0,\dots,c_3$, $\gamma$ are  real parameters  and  $\varepsilon$ characterizes the dispersion.

There is known (see e.g. \cite{ELY}) that the family (\ref{1}) contains only three special cases that satisfy the asymptotic integrability condition: the Korteweg-de Vries, the Camassa-Holm, and the Degasperis-Procesi    equations. More in detail:

1. Obviously, if we set $\alpha=c_1=c_2=0$ then we obtain the KdV equation
\begin{equation}
\frac{\pa u}{\pa t}+\frac{\pa}{\pa x}\left\{c_0u+c_1u^2
+\gamma\ve^2\frac{\pa^2 u}{\pa x^2}\right\}=0,\label{2}
\end{equation}
which describes   the wave propagation at the free surface of shallow water under the influence of gravity.

2. For $c_1=3c_3/(2\alpha^2)$, $c_2=c_3/2$, $\gamma=0$, and $v=c_3u$ Eq. (\ref{1}) becomes the Camassa-Holm equation
 modeling the propagation of shallow water waves over a flat bottom,
\begin{equation}
\frac{\pa }{\pa t}\left\{v-\alpha^2\ve^2\frac{\pa^2 v}{\pa x^2}\right\}+\frac{\pa}{\pa x}\left\{c_0v+\frac{3}{2\alpha^2}v^2
-\ve^2\left(\frac12\Big(\frac{\pa v}{\pa x}\Big)^2+v\frac{\pa^2 v}{\pa x^2}\right)\right\}=0.\label{3}
\end{equation}

3. If $c_1=2c_3/\alpha^2$, $c_2=c_3$, $c_0=\gamma=0$, and $v=c_3u$, then  Eq. (\ref{1}) becomes the Degasperis-Procesi equation,
\begin{equation}
\frac{\pa }{\pa t}\left\{v-\alpha^2\ve^2\frac{\pa^2 v}{\pa x^2}\right\}+\frac{\pa}{\pa x}\left\{\frac{2}{\alpha^2}v^2
-\ve^2\left(\Big(\frac{\pa v}{\pa x}\Big)^2+v\frac{\pa^2 v}{\pa x^2}\right)\right\}=0.\label{4}
\end{equation}
It is known that the KdV equation and the Camassa-Holm equation for $c_0>0$ admit smooth solitary wave solutions called "solitons", see
Fig.1.
\begin{figure}[H]
\centering
\includegraphics[height=2in,width=3.5in]{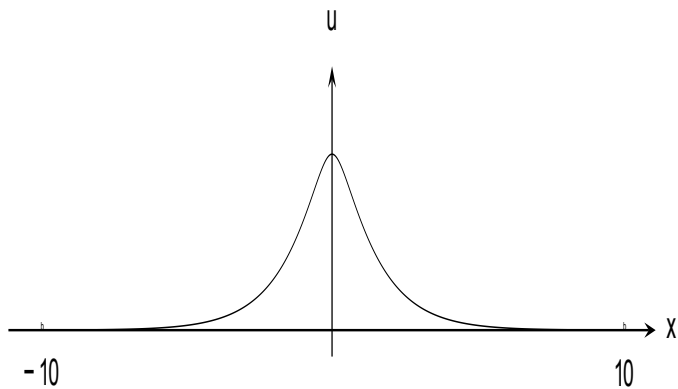}
\caption{Soliton solution of the  Camassa-Holm equation for $c_0>0$}
\label{f2}
\end{figure}
Conversely,
the Degasperis-Procesi equation and the Camassa-Holm equation for $c_0=0$ have continuous solitary wave solutions called "peacons", see Fig.2.
\begin{figure}[H]
\centering
\includegraphics[width=13cm]{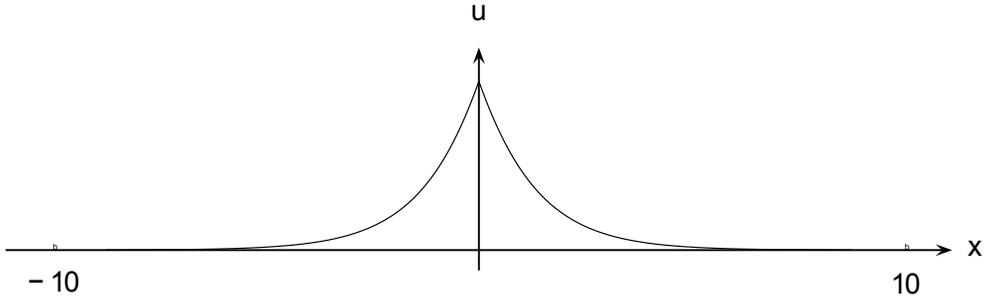}
\caption{Peacon  solution of the Camassa-Holm equation for $c_0=0$}
\label{f3}
\end{figure}
Moreover, it is known that the solitary wave solutions of the equations (\ref{3}) and (\ref{4}) interact elastically, that is in the same manner as the KdV solitons.

However, the special cases (\ref{2}) - (\ref{4}) exhaust that's all what is known about the family (\ref{1}). To begin the study of the wave propagation for nonintegrable versions of (\ref{1})  we should separate firstly two basic situations: smooth and non-smooth traveling solutions. We will do it in Section 2. The next question is the scenario of the solitary wave interaction. We shall discuss it in Section 3 for the case of solitons.

\section{Solitary wave solution}
Let us set the ansatz
\begin{equation}\label{5}
u=A\om\big(\beta(x-Vt)/\varepsilon\big),
\end{equation}
where $\om(\eta)$ is a smooth even function such that
\begin{align}
&\om(\eta)\to0\quad\text{as}\quad \eta\to\pm\infty,\label{6}\\
&\om(0)=1,\label{7}
\end{align}
the amplitude $A>0$ is a  free parameter, and the parameters
$\beta$, $V$ should be determined.  In what follows we assume that
\be\label{8}
\gamma\geq0,\quad c_0\geq0,\quad\alpha>0,\quad c_k>0,\quad k=1,2,3.
\ee
Substituting  (\ref{5}) into Eq. (\ref{5}), integrating, and using (\ref{6}), we obtain the second order ODE
\begin{align}
\left\{1-\frac{c_2c_4A }{\gamma+\alpha^2V}\om\right\}\frac{d^2 \om}{d \eta^2}&=\frac{c_2A }{\gamma+\alpha^2V}\left(\frac{d \om}{d \eta}\right)^2\notag\\
&+\frac{1}{\beta^2(\gamma+\alpha^2V)}\big((V-c_0)\om-c_1A\om^2\big).\label{9}
\end{align}
Next we define the auxiliary parameter $\beta$,
\begin{equation}\label{10}
\beta^2=c_1/c_3,
\end{equation}
rescaling the function $\om$,
\begin{equation}\label{11}
 W=p\om,\quad p=c_3A/(\gamma+\alpha^2V),
\end{equation}
and denote constants
\begin{equation}
r=c_3/(c_2+c_3), \quad q=c_3(V-c_0)/(\gamma+\alpha^2V).\label{12}
\end{equation}
Using (\ref{10}) - (\ref{12}) we deduce that $W$ satisfies the equation
\begin{equation}\label{13}
(1-W)\frac{d^2 W}{d\eta^2}=\frac{1-r}{r}\left(\frac{d W}{d\eta}\right)^2+qW-W^2.
\end{equation}
The next step is the substitution
\begin{equation}\label{14}
 W(\eta)=1-g(\eta)^r,
\end{equation}
which allows us to eliminate the first derivatives from the model equation (\ref{13}). Taking into account the condition (\ref{7}) and the property of being even, $g(-\eta)=g(\eta)$, we pass to the "boundary" problem
\begin{align}
&r\frac{d^2 g}{d\eta^2}=g-(2-q)g^{1-r}+(1-q)g^{1-2r},\quad \eta\in(0,\infty),\label{15}\\
& g^r\big|_{\eta=0}=1-p,\quad g|_{\eta\to\infty}=1.\label{16}
\end{align}
Notice that the correctness of (\ref{16}) implies the assumption
\begin{equation}
\text{if}\quad p>1, \quad \text{then}\quad r=(2k+1)/(2l+1),\quad k,l\in \mathbb{Z}.\label{17}
 \end{equation}
 Now we integrate (\ref{15}) and pass to the first order ODE
\begin{equation}\label{18}
 r\left(\frac{d g}{d\eta}\right)^2=F(g),\quad \eta\in(0,\infty);\quad
 g|_{\eta=0}=g_*,
\end{equation}
where
\begin{align}
&F(g)=g^2-2\frac{2-q}{2-r}g^{2-r}+\frac{1-q}{1-r}g^{2-2r}-C,\label{19}\\
&C=1-2\frac{2-q}{2-r}+\frac{1-q}{1-r},\quad g_*=(1-p)^{1/r}.\label{20}
\end{align}
Considering $\eta>>1$ we write $g=1-w$ and obtain from (\ref{18})-(\ref{20})
$$
\left(\frac{d w}{d\eta}\right)^2=q\,w^2.
$$
Thus
$$
g\to 1-e^{-\sqrt{q}\eta}\quad\text{as}\quad \eta\to\infty.
$$
Therefore, the solution of the problem (\ref{18}) exists, is unique for $g>0$, and satisfies the conditions (\ref{16}).

We now consider the even continuation $\widetilde{g}(\eta)$ of $g$ for negative $\eta$. Obviously, $\widetilde{g}\in C^\infty(\mathbb{R})$ if and only if
\begin{equation}\label{21}
\frac{d g}{d\eta}\Big|_{\eta=0}=0.
\end{equation}
Furthermore, since $F(1)=dF/dg|_{g=1}=0$ and $d^2F/dg^2|_{g=1}>0$, the equation
\begin{equation}\label{22}
F(g)=0
\end{equation}
has a solution $g_*\in(0,1)$ if and only if $C>0$. The last inequality is equivalent to the following assumption:
\begin{equation}\label{23}
r>q.
\end{equation}
On the other hand the initial condition in (\ref{18}) implies the relation
\begin{equation}\label{24}
V=\frac{1}{\alpha^2}\left(\frac{c_3}{1-g_*^r}A-\gamma\right).
\end{equation}
This allows to rewrite the coefficient $q$ in (\ref{19}) as a function on $A$ and the parameters $\alpha$, $c_0,\dots,c_3$, $\gamma$; therefore to find the solution $g_*$ of the equation (\ref{22}) as a function on $A$ and the parameters of the model (\ref{1}). Representing (\ref{23}) in the explicit form we obtain the conclusion
\begin{teo}
Under the assumption (\ref{8}) we assume that
\begin{equation}\label{25}
c_3-A^{-1}(\gamma+c_0\alpha^2)(1-g_*^r)<c_1r\alpha^2.
\end{equation}
Then the equation (\ref{9}) has the unique smooth solution what is even and satisfies the conditions  (\ref{6}), (\ref{7}).
\end{teo}
\textbf{Example 1.} For the Camassa-Holm equation (\ref{3}) $r=2/3$ and (\ref{22}) is a cubic equation. Thus
$$
g_*= \big(1+c_3A/(c_0\alpha^2)\big)^{-3/2}\quad\text{if}\quad c_0>0\quad \text{and}\quad g_*=0\quad\text{if}\quad c_0=0.
$$
Respectively the condition (\ref{25}) is satisfied for $c_0>0$ and it is broken for $c_0=0$. In the last case
$$
\frac{d \om}{d\eta}\Big|_{\eta=0}=-\sqrt{2(1-q)}/p\neq 0,
$$
therefore $\om(\eta)$ is a continuous function only. Fig.3 depicts the $F(g)$ graph in the case $c_0=1$, $A=2$,  $c_3=2$ and $\alpha=1$.
\begin{figure}[H]
\centering
\includegraphics[width=9cm]{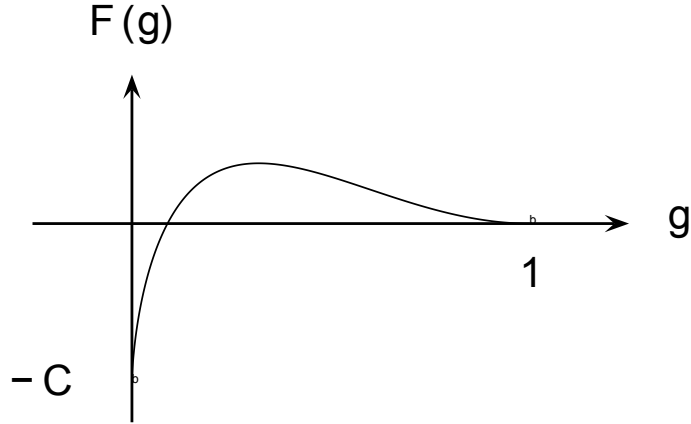}
\caption{Behavior of the function $F(g)$ for the Camassa-Holm equation with $c_0=1$}
\label{f1}
\end{figure}

\textbf{Example 2.} For the Degasperis-Procesi equation (\ref{4}) the condition (\ref{25}) is violated and
$$
\frac{d \om}{d\eta}\Big|_{\eta=0}=-\sqrt{(1-q)}/p\neq 0.
$$
Fig.2 demonstrates the graph of the peacon $\om(\eta)$ for this equation.

\textbf{Example 3}
Now let $c_0=\gamma=0$ and $\alpha^2c_1>c_2+c_3$. Then $q=c_3/\alpha^2c_1<r$ and $g_*$ doesn't depend on $V$. Thus
\begin{equation}
V=c_3A/\{(1-{g_*}^r)\alpha^2\}.\label{26}
\end{equation}

\textbf{Example 4}
Let $c_3=4c_2$. Setting $z=g^r$, $r={2/5}$, we transform the equation (\ref{22}) to the form
\begin{equation}
F=(1-z)^2f=0,\quad f=z^3+2z^2-\frac13(1-5q)z-\frac45(4-5q).\label{27}
\end{equation}
Solving the cubic equation $f=0$ we find the root $z_*=z_*(V)$. This and  (\ref{24}) imply the equality
\begin{equation} \label{27a}
A=\mathfrak{A}(V),\quad \mathfrak{A}=c_3^{-1}(\gamma+\alpha^2V)\big(1-z_*(V)\big).
\end{equation}
Simple calculations show that  $d\mathfrak{A}/dV|_{q=0}>0$. Thus,   (\ref{27a}) allows us to define the velocity as a function of  the amplitude at least for $V-c_0<<1$.

Similar result can be obtained in the case $c_2=3c_3/2$.
\section{Two-soliton asymptotic solution}
Obviously, there is not any hope to find both the exact multi-soliton solution to (\ref{1}) and an
asymptotics in the classical sense. So, we will treat $\ve$ as a small parameter and construct a
weak asymptotic solution.
The Weak Asymptotics Method (see e.g. \cite{DanShel1} - \cite{Colombeau} and references therein) takes
into account the fact that soliton-type solutions which are smooth for $\ve>0$ become non-smooth in the limit as $\ve\to0$. Thus, it is
possible to treat such solutions  as a mapping $\mathcal{C}^\infty
(0, T; \mathcal{C}^\infty (\mathbb{R}_x^1))$ for $\ve=\const>0$ and only
as $\mathcal{C} (0, T; \mathcal{D}' (\mathbb{R}_x^1))$ uniformly
in $\ve\ge0$. Accordingly, the remainder should be small in the
weak sense. The main advantage of the method is such that we can ignore the real shape of the colliding waves but look for (and find) exceptionally their main characteristics. For the equation (\ref{1}) they are the amplitudes and trajectories of the waves.

Originally, such idea had been suggested by Danilov\&Shelkovich
for shock wave type solutions (\cite{DanShel1}, 1997), and by  Danilov\&Omel'yanov for soliton type solutions (\cite{DanOm}, 2003). Later the method has
been developed and adapted   for many other problems (V. Danilov, G. Omel'yanov, V. Shelkovich, D. Mitrovic, M. Colombeau and others, see e.g. \cite{DanOmShel} - \cite{Colombeau}).

Notice finally that the treatment (Omel'yanov \cite{Om2, Om3}, 2015) of weak asymptotics as functions
which satisfy some conservation or balance laws  takes us back to the ancient Whitham's idea to construct one-phase
asymptotic solution satisfying a Lagrangian. Now, for essentially nonintegrable equations and multi-soliton solutions, we
use the appropriate number of the laws and satisfy them in the weak sense.

Let us apply these ideas for  the problem of two soliton interaction in the Degasperis-Procesi  model (\ref{1}).
We set the initial data
\begin{equation} \label{28}
u|_{t=0} = \sum_{i=1}^2 A_i\omega \big(\beta (x - x_i^0)/\varepsilon\big).
\end{equation}
Here
$A_2>A_1>0$, $x_1^0 >x_2^0$,  $\beta=\sqrt{c_1/c_3}$,  and we assume that  the trajectories $x = V_i t + x_i^0$
have a joint point $x = x^*$ at a time instant $t = t^*$, where $V_i$ are defined in the same manner as in
 (\ref{24}).

To construct the weak asymptotic solution we start with the following definition of the smallness in the weak sense \cite{DanOm, Om3}:

{\bf Definition 1}

{\it A function}  $v(t, x, \varepsilon)$ {\it is said to be of the value} $
O_{\mathcal{D}'}(\varepsilon^\varkappa)$ {\it if the relation}
$$ \int_{-\infty}^\infty v(t, x, \varepsilon )\psi(x) dx = O(\varepsilon^\varkappa)$$
{\it holds uniformly in} $t$ {\it for any test function} $\psi \in \mathcal{D}
(\mathbb{R}_x^1)$. {\it The right-hand side here is a} $\mathcal{C}^\infty${\it-function for}
$\varepsilon=\const > 0$ {\it and  a piecewise continuous function
uniformly in} $\varepsilon \geq 0$.

Next we write two  associated with (\ref{1}) conservation and balance laws  in the differential form:
\begin{equation}\label{29}
\frac{\partial Q_j}{\partial t}+\frac{\partial P_j}{\partial x}+\ve^{-1}K_j=O_{\mathcal{D}'}(\varepsilon^2),\quad j=1,2,
\end{equation}
where
\begin{align}
&Q_1=u,\quad P_1=c_0u+c_1u^2-(c_2-c_3)(\ve u_x)^2,\quad K_1=0,\label{30}\\
&Q_2=u^2+\alpha^2(\ve u_x)^2,\quad P_2=\mathbb{P}_2+2\alpha^2\ve^2 u_x u_t,\quad
K_2=(2c_2-c_3)(\varepsilon u_x)^3,\label{31}\\
&\mathbb{P}_2=c_0u^2+\frac43c_1u^3-\big(3\gamma +(2c_2-5c_3)u\big)(\ve u_x)^2,\label{32}
\end{align}
where subscripts denote partial derivatives.

Following  \cite{DanOm, Om3}, we define two-soliton weak asymptotics:

{\bf Definition 2}

{\it A sequence} $u(t, x, \varepsilon )$, {\it belonging to}
$\mathcal{C}^\infty (0, T; \mathcal{C}^\infty (\mathbb{R}_x^1))$
{\it for} $\varepsilon =\const> 0$ {\it and belonging to} $\mathcal{C} (0, T;
\mathcal{D}' (\mathbb{R}_x^1))$ {\it uniformly in} $\varepsilon\geq0$, {\it is
called a weak asymptotic} mod $ O_{\mathcal{D}'}(\varepsilon^2)$
{\it solution of (\ref{1}) if the relations (\ref{29}) hold uniformly in}
$t$  {\it with the accuracy} $O_{\mathcal{D}'}(\varepsilon^2)$.

Next we present the ansatz
 as the sum of two distorted solitons, that is:
\begin{equation}\label{33}
u=\sum_{i=1}^2G_i\om\big(\beta(x-\vp_i)/\varepsilon\big),
\end{equation}
where
\be
G_i=A_i+S_i(\tau),\,\vp_i=\vp_{i0}(t)+\ve\vp_{i1}(\tau),\,\tau=\beta_1\big(\vp_{20}(t)-\vp_{10}(t)\big)/\ve,\label{34}
\ee
$\vp_{i0}=V_it+x_{i0}$ describe the trajectories of the non-interacting  waves
(\ref{5}) with the amplitudes $A_i$.
Next  we suppose  that $S_i(\tau)$, $\vp_{i1}(\tau)$ are smooth functions such that
\begin{align}
&S_i\to0 \quad \text{as}\quad \tau\to\pm\infty,\label{35}\\
&\vp_{i1}\to0 \quad \text{as}\quad \tau\to-\infty,\quad
\vp_{i1}\to \vp_{i1}^\infty=\const_i \quad \text{as}\quad
\tau\to+\infty.\label{36}
\end{align}
It is obvious that the existence of the weak asymptotics  (\ref{33}) with the properties (\ref{35}), (\ref{36})
implies that the solitary waves  interact like the KdV solitons at least in the leading term.

To construct the asymptotics we should calculate  the weak expansions of the terms from the left-hand sides of the relations (\ref{29}).
It is easy to check that
\be
u=\ve\beta^{-1}a_{1} \sum_{i=1}^2G_i\delta(x-\vp_i)+O_{\mathcal{D}'}(\ve^3),\label{37}
\ee
where $\delta(x)$ is the Dirac delta-function. Here and in what follows we use the notation
\begin{equation}\label{38}
a_{k}\ipd\int_{-\infty}^\infty\big(\om(\eta)\big)^k d\eta, \quad k>0,\qquad
a'_{2}\ipd\int_{-\infty}^\infty\big(\om'(\eta)\big)^2 d\eta.
\end{equation}
At the same time for any even $F(u,\ve u_x)\in C^1$, $F\big(u(-x),\ve u_x(-x)\big)=F\big(u(x),\ve u_x(x)\big)$, we have
\begin{align}
&\int_{-\infty}^\infty F\left(\sum_{i=1}^2 G_i\omega \left( \beta \frac{x - \vp_i
}{\varepsilon}\right), \beta\sum_{j=1}^2 G_j\omega' \left( \beta \frac{x - \vp_j
}{\varepsilon}\right)\right)\psi(x)dx\notag\\
&=\frac{\ve}{\beta}\int_{-\infty}^\infty \sum_{i=1}^2F\big(A_i\om(\eta),\beta A_i\om'(\eta)\big)
\psi(\vp_i+\ve\frac\eta\beta)d\eta\label{39}\\
&+\frac{\ve}{\beta}\int_{-\infty}^\infty\Big\{ F\Big(\sum_{i=1}^2G_i\om(\eta_{i2}),\beta\sum_{j=1}^2G_j\om'(\eta_{j2})\Big)\notag\\
&-\sum_{i=1}^2F\big(A_i\om(\eta_{i2}),\beta A_i\om'(\eta)\big) \Big\}\psi(\vp_2+\ve\frac\eta\beta)d\eta,\notag
\end{align}
where
\be
\eta_{12}=\eta-\sigma, \quad \eta_{22}=\eta,\quad\sigma=\beta(\vp_1-\vp_2))/\ve.\label{40}
\ee
We take into account that the second integrand in right-hand side (\ref{40}) vanishes exponentially fast as $|\vp_1-\vp_2|$  grows,
thus, its main contribution is at the point $x^*$.  We write
\be
\vp_{i0}=x^*+V_i(t-t^*)=x^*+\ve \frac{V_i}{\dpsi}\tau \quad \text{and} \quad \vp_{i}=x^*+\ve\chi_{i},\label{41}
\ee
where $\dpsi=\beta(V_2-V_1)$, $\chi_i=V_i\tau/\dpsi + \vp_{i1}$.
It remains to apply the formula
\be
f(\tau)\delta(x-\vp_i)=f(\tau)\delta(x-x^*)-\ve\chi_if(\tau)\delta'(x-x^*)+O_{\mathcal{D}'}(\ve^2),\label{42}
\ee
which holds for each $\vp_i$ of the form (\ref{41}) with slowly increasing $\chi_i$ and for $f(\tau)$ from the Schwartz space.
 Moreover, the second term in (\ref{42}) is $O_{\mathcal{D}'}(\ve)$. Thus, under the assumptions (\ref{35}), (\ref{36})
  we obtain the weak asymptotic expansion of $F(u)$ in the final form:
\be
F(u)=\frac{\ve}{\beta} \sum_{i=1}^2a_{F,i}\delta(x-\vp_{i})+\frac{\ve}{\beta}\mathfrak{R}_{F}^{(0)}\delta(x-x^*)
+O_{\mathcal{D}'}(\varepsilon^2),\label{43}
\ee
where
\begin{align}
 a_{F,i}=\int_{-\infty}^\infty &F\big(A_i\om(\eta),\beta A_i\om'(\eta)\big)d\eta, \label{44}\\
\mathfrak{R}_{F}^{(n)}=
\int_{-\infty}^\infty\eta^n\Big\{& F\Big(\sum_{i=1}^2G_i\om(\eta_{i2}),\beta\sum_{j=1}^2G_j\om'(\eta_{j2})\Big)\label{45}\\
&-\sum_{i=1}^2F\big(A_i\om(\eta_{i2}),\beta A_i\om'(\eta_{i2})\big) \Big\}d\eta,\quad n=1,2.\notag
\end{align}
Note that to define $\pa Q_2/\pa t \mod O_{\mathcal{D}'}(\varepsilon^2)$ it is necessary to
calculate $Q_2$ with the precision $O_{\mathcal{D}'}(\varepsilon^3)$. Thus, transforming (\ref{37}) with the help
of (\ref{42}) and using (\ref{43}) with $F(u)=Q_2$, we obtain modulo $O_{\mathcal{D}'}(\ve^3)$:
\begin{align}
u=a_{1}\frac{\ve}{\beta} \sum_{i=1}^2A_{i}\delta(x-\vp_i)
&+a_{1}\frac{\ve}{\beta} \sum_{i=1}^2S_{i}\Big\{\delta(x-x^*)-\ve\chi_i\delta'(x-x^*)\Big\},\label{46}\\
Q_2=\frac{\ve}{\beta}\sum_{i=1}^2a_{Q_2,i}\delta(x-\vp_i)
&+\frac{\ve}{\beta} \mathfrak{R}_{Q_2}^{(0)}
\delta(x-x^*)\label{47}\\
&-\frac{\ve^2}{\beta}\big\{\chi_2\mathfrak{R}_{Q_2}^{(0)}+\beta^{-1}\mathfrak{R}_{Q_2}^{(1)} \big\}\delta'(x-x^*).\notag
\end{align}
In the same manner we derive
\begin{align}
&\ve^2u_xu_t=-\ve a'_{2}\beta\sum_{i=1}^2V_iA^2_{i}\delta(x-\vp_i)
-\ve a'_{2}\mathfrak{L}\delta(x-x^*),\label{48}\\
&(\ve u_x)^3=\frac{\ve}{\beta} \mathfrak{R}_{K_2}^{(0)}\delta(x-x^*)-\ve^2 \beta a^{(1)}_{3}\beta\sum_{i=1}^2A^3_{i}\delta'(x-\vp_i)\notag\\
&-\frac{\ve^2}{\beta} \left\{\chi_2\mathfrak{R}_{K_2}^{(0)}+\beta^{-1}\mathfrak{R}_{K_2}^{(1)}\right\}\delta'(x-x^*),\label{49}
\end{align}
where
\begin{align}
&\lambda_{(k,l)}=\frac{1}{a'_2} \int_{-\infty}^\infty \om^{(k)}(\eta_{12})\om^{(l)}(\eta)d\eta,\quad
a^{(1)}_{3}=\int_{-\infty}^\infty\eta\big(\om'(\eta)\big)^3 d\eta.\label{50}\\
&\mathfrak{L}=\dpsi\beta\sum_{i=1}^2\frac{d\vp_{i1}}{d\tau}(G_i^2-A_i^2)
-\dpsi\left(G_{1}\frac{dS_{2}}{d\tau}-G_{2}\frac{dS_{1}}{d\tau}\right)
\lambda_{(1,0)}\notag\\
&+\beta G_1G_2(\dot{\vp_1}+\dot{\vp_1})\lambda_{(1,1)}.
\end{align}
Substituting (\ref{43})-(\ref{49}) into
(\ref{29}) we obtain linear combinations of $\delta(x-x^*)$,
$\ve\delta'(x-\vp_{i})$, $i=1,2$,  and
$\ve\delta'(x-x^*)$ (see also \cite{DanOm, DanOmShel, Om3}). Therefore, we pass to the following system of equations:
\begin{align}
&a_{1}V_iA_{i}-a_{P_1,i}=0, \, a_{2}V_i\beta^2A_{i}^{2}+a_{Q_2,i}V_i-a_{\mathbb{P}_2,i}+a^{(1)}_{3}\beta_i^2A_{i}^3=0, \, i=1,2,\label{51}\\
&\sum_{i=1}^2S_{i}=0,\quad \dpsi \frac{d}{d\tau}\mathfrak{R}_{Q_2}^{(0)}+\mathfrak{R}_{K_2}^{(0)}=0,\label{52}\\
&a_1\dpsi \frac{d}{d\tau}\sum_{i=1}^2\Big\{A_{i}\vp_{i1}+\chi_iS_i\Big\}=f,\label{53}\\
&\dpsi\frac{d}{d\tau}\Big\{\sum_{i=1}^2a_{Q_2,i}\vp_{i1}+\chi_2\mathfrak{R}_{Q_2}^{(0)}+\beta^{-1}\mathfrak{R}_{Q_2}^{(1)}\Big)=F,\label{54}
\end{align}
where
\be
f=\mathfrak{R}^{(0)}_{P_1}, \quad  F=\mathfrak{R}^{(0)}_{P_2}-a'_{2}\mathfrak{L}-\chi_2\mathfrak{R}^{(0)}_{K_2}-\beta^{-1}\mathfrak{R}^{(1)}_{K_2}.\label{55}
\ee
An analysis of (\ref{51}) implies the following statement:
\begin{lem}
The algebraic equations (\ref{51})  with $\beta=\sqrt{c_1/c_3}$ imply again the  relation
(\ref{24}) between $A_i$ and $V_i$.
\end{lem}
As for  (\ref{52})-(\ref{54}), this system should be investigated in detail. Now we can formulate a previous result only
\begin{teo} \label{teo1}
Let the assumptions (\ref{8}), (\ref{25}) be satisfied.
Presuppose also that the equations (\ref{52})-(\ref{54}) admit a solution with the properties (\ref{3}), (\ref{36}). Then the solitary wave collision  in the problem
(\ref{1}), (\ref{28}) preserves the elastic scenario with accuracy
$O_{\mathcal{D}'} (\varepsilon^2)$ in the sense of Definition 2.
\end{teo}

\end{document}